\documentstyle[twoside,11pt,leqno]{article}

   \def\H{{\cal H}} 
\def\D{{\texttt{D}}}
 \def\a{\rm{asc}}
 \def\N{{\texttt{N}}}
\def\C{\texttt{C}} 

\def\P{{\rm {H_c}}}

\def\B{B({\cal H})} 
\def\asc{ \textrm{asc}} 

\newtheorem{df}{Definition}[section]
\newtheorem{thm}[df]{Theorem} \newtheorem{pro}[df]{Proposition}
\newtheorem{cor}[df]{Corollary} 
\newtheorem{rema}[df] {Remark} 
\def\sfstp{{\hskip-1em}{\bf.}{\hskip1em}}

\def\subject#1{\renewcommand{\thefootnote}{}\footnote
{AMS(MOS) subject classification (2010). Primary: {#1}}}

\def\keywords#1{\renewcommand{\thefootnote}{}\footnote
{Keywords: {#1}}}

\def\enddemo{\qed \endtrivlist} \expandafter\let\csname
enddemo*\endcsname=\enddemo

\def\qedsymbol{\ifmmode\bgroup\else$\bgroup\aftergroup$\fi
\vcenter{\hrule\hbox{\vrule
height.5em\kern.5em\vrule}\hrule}\egroup}
\def\qed{\ifmmode\else\unskip\nobreak\fi\quad\qedsymbol}

\title{\bf  Power bounded $m$-left invertible operators}
\author{\normalsize B.P. Duggal and C.S. Kubrusly}
\date{}

\begin{document}

\maketitle \thispagestyle{empty} \vskip-16pt

\subject{Primary47A05, 47A55; Secondary 47A80, 47A10.} \keywords{Hilbert space, Power bounded operator, $m$-left invertible operator, $m$-isometric operator, $(m,C)$-isometric operator, similar to an isometry.}\footnote{}
\begin{abstract} A Hilbert space operator $S\in\B$ is left $m$-invertible by $T\in\B$ if $$\sum_{j=0}^m{(-1)^{m-j}\left(\begin{array}{clcr}m\\j\end{array}\right)T^jS^j}=0,$$  $S$ is $m$-isometric if $$\sum_{j=0}^m{(-1)^{m-j}\left(\begin{array}{clcr}m\\j\end{array}\right){S^*}^jS^j}=0$$ and $S$ is $(m,C)$-isometric for some conjugation $C$ of $\H$ if $$\sum_{j=0}^m{(-1)^{m-j}\left(\begin{array}{clcr}m\\j\end{array}\right){S^*}^jCS^jC}=0.$$ If a power bounded operator $S$ is left invertible by a power bounded operator $T$, then $S$ (also, $T^*$) is similar to an isometry. Translated to $m$-isometric and $(m,C)$-isometric operators $S$ this implies that $S$ is $1$-isometric, equivalently isometric, and (respectively) $(1,C)$-isometric.
\end{abstract}


\section {\sfstp Introduction}Given a complex infinite dimensional Hilbert space $\H$ , let $\B$ denote the algebra of bounded linear transformations, equivalently operators,
on $\H$  into itself. Given operators $S, T\in\B$, let $$ P_m(S,T)=\sum_{j=0}^m{(-1)^{m-j}\left(\begin{array}{clcr}m\\j\end{array}\right)T^jS^j}.$$  We say that $T$ is a left $m$-inverse of $S$ (equivalently, $S$ is left $m$-invertible by $T$) 
for some integer $m>0$ if $$P_m(S,T)=0$$ \cite{{OAMS},{DM},{CG}}. 
Left $m$-invertible operators occur quite naturally, and the class of  $m$-isometric operators, i.e. operators $S\in \B$ such that
$$P_m(S,S^*)=\sum_{j=0}^m{(-1)^{m-j}\left(\begin{array}{clcr}m\\j\end{array}\right){S^*}^jS^j}=0,$$
of Agler and Stankus \cite{AS} is an important widely studied example of operators left $m$-invertible by their adjoint. A generalization of $m$-isometric operators, which has been considered in the recent past \cite{CKL}, is that of $(m,C)$-isometric 
operators. Here an operator $S\in\B$ is $(m,C)$-isometric for some conjugation $C$ of $\H$ if
$$P_m(CSC,S^*)=\sum_{j=0}^m{(-1)^{m-j}\left(\begin{array}{clcr}m\\j\end{array}\right)S^{*j}CS^jC}=0.$$ 
(Recall that a conjugation $C$ of $\H$ is an antilinear operator such that $C^2=I$ and $<Cx,Cy>=<y,x>$ for all $x,y\in\H$.)

\

An operator $S\in\B$ is power bounded if there exists a scalar $M>0$ such that $$\sup_{n\in\N}{||S^n||< M}.$$ It is immediate from the definition that if $S\in\B$ is power bounded, then the spectral radius $$r(S)=\lim_{n\rightarrow\infty}{||S^n||^{\frac{1}{n}}\leq 1}$$ and the spectrum $\sigma(S)$ of $S$ satisfies $\sigma(S)\subseteq\D$ ($=\{\lambda\in\C: |\lambda|\leq 1\}$). Recall from \cite[Theorem 2.4]{D} that an $m$-isometric operator is power bounded if and only if it is isometric.

\

Given a positive operator $A\in\B$, $A\geq 0$, let $||.||_A$ denote the semi-norm
$$||x||_A^2=<x,x>_A=<Ax,x>, \hspace{2mm} x\in\H.$$
(Then $||.||_A$ is a norm on $\H$ if and only if $A$ is injective.) An operator $S\in\B$ is said to be $A$-isometric if $S^*AS=A$ \cite{ACG} and $S$ is  an $(A,m)$-isometry \cite{D} if $$P_m(A;S,S^*)=\sum_{j=0}^m{(-1)^{m-j}\left(\begin{array}{clcr}m\\j\end{array}\right){S^{*j}A S^j}}=0.$$ This paper considers left $m$-invertible operators such that both the operator $S\in\B$ and its left $m$-inverse $T\in\B$ are power bounded. It is proved that there exist positive invertible operators $P_i$, $P_i>0$ ($i=1,2$), such that $S=P_1V_1P_1^{-1}$ and $T^*=P_2V_2P_2^{-1}$ for some isometries $V_i$, $i=1,2$. Translated to $m$-isometric  and $(m,C)$-isometric $S$ this means that: {\em a power bounded $m$-isometric operator is isometric and a power bounded $(m,C)$-isometric operator is $(1,C)$-isometric}.

\section {\sfstp Results} Given an operator $A\in\B$, we write $A-\lambda$ for $A-\lambda I$, $\lambda\in\C$.  $A$ has SVEP, {\em the single-valued
	extension property}, at $\lambda_0\in\C$ if for every open
disc ${\D}_{\lambda_0}$ centered at $\lambda_0$ the only
analytic function $f:{\D}_{\lambda_0}\longrightarrow \H$
satisfying $(A-\lambda)f(\lambda)=0$ is the function $f\equiv 0$
\cite{{A},{LN}} . Every operator $A$ has SVEP at points in its resolvent set
$\rho(A)=\C\setminus \sigma(A)$ and  on the boundary $\partial\sigma(A)$
of the spectrum $\sigma(A)$. {\em We say that $A$ has SVEP on a set
	$\Xi$ if it has SVEP at every $\lambda\in \Xi$.} The {\em ascent of
	$A$}, $\asc(A)$, is the
least non-negative integer $n$ such that $A^{-n}(0)=A^{-(n+1)}(0)$: If no such integer exists, then
$\asc(A)=\infty$. It is well known  that
$\asc(A)<\infty$ implies $A$ has SVEP at $0$ \cite{{A},{LN}}.

\

For $A. B\in\B$, let $\triangle_{A,B}\in B(\B)$ denote the elementary operator $\triangle_{A,B}(X)=AXB-X=(L_AR_B-I)(X)$, where $L_A$ and $R_B\in B(\B)$ are the operators $L_A(X)=AX$ and $R_B(X)=XB$ (of left multiplication by $A$ and, respectively, right multiplication by $B$). It is known, see for example \cite{Sh}, that if $A,B$ are normal operators, then ($\triangle_{A,B}^{-1}(0)\subseteq\triangle^{-1}_{A^*,B^*}(0)$, consequently) $\a(\triangle_{A,B})\leq 1$. 

\

$A\in\B$ is a $C_{0.}$,  respectively $C_{1.}$, operator if $$ \lim_{n\rightarrow\infty}{||A^nx||=0}\hspace{2mm}\mbox{for all}\hspace{2mm}x\in\H,$$ $$\mbox{respectively,} \inf_{n\in\N}{||A^nx||>0}\hspace{2mm}\mbox{for all}\hspace{2mm}0\neq x\in\H;$$
 $A\in C_{.0}$ (resp., $A\in C_{.1}$) if $A^*\in C_{0.}$ (resp., $A^*\in C_{1.}$) and $A\in C_{\alpha\beta}$ if $A\in C_{\alpha .}\cap C_{. \beta}$ ($\alpha,\beta=0,1$). It is well known \cite{K} that every power bounded operator $A\in\B$ has an upper triangular matrix representation $$A=\left(\begin{array}{clcr}A_1 & A_0\\0 & A_2\end{array}\right)\in B(\H_1\oplus\H_2)$$ for some decomposition $\H=\H_1\oplus\H_2$ of $\H$ such that $A_1\in C_{0.}$ and $A_2\in C_{1.}$. Recall that every isometry $V\in\B$ has a direct sum decomposition $$V=V_c\oplus V_u\in B(\H_c\oplus\H_u), V_c\in C_{10}\hspace{2mm}\mbox{and}\hspace{2mm} V_u\in C_{11}$$ into its completely non-unitary (i.e., unilateral shift) and unitary parts.
 
 \
 Let $\delta_{A,B}\in B(\B)$ denote the generalized derivation $\delta_{A,B}(X)=AX-XB$. Recall from \cite{DK} that $A\in\B$ satisfies (the Putnam-Fuglede) property $\rm{PF}(\triangle)$  (resp., $\rm{PF}(\delta)$),  $A\in\rm{PF}(\triangle)$ (resp., $A\in\rm{PF}(\delta)$), if (either $A$ is trivially unitary, or) for every isometry $V\in\B$ for which $\triangle_{A,V^*}(X)=0$ (resp. $\delta_{A^*,V}(X)=0$) has a non-trivial solution $X\in\B$, $X$ is also a solution of $\triangle_{A^*,V}(X)=0$ (resp., $\delta_{A^*,V}(X)=0$). The following theorem is \cite[Corollary 2.5]{DK} (see also \cite{P}). Let $d_{A,B}$ denote either of $\triangle_{A,B}$ and $\delta_{A,B}$, and, correspondingly, let $\rm{PF}(d)$ denote either of  $\rm{PF}(\triangle)$ and $\rm{PF}(\delta)$. (Recall from \cite[Theorem 2.1]{DK} that $A\in \rm{PF}(\triangle)$ if and only if $A\in \rm{PF}(\delta)$.)
\begin{thm}
	\label{thm01} A power bounded operator satisfies property $\rm{PF}(d)$ if and only if it is the direct sum of a unitary with a $C_{.0}$ operator.
\end{thm} The $\rm{PF}(d)$ property implies range-kernel orthogonality (i.e., if $d_{A,V^*}^{-1}(0)\subseteq d_{A^*,V}^{-1}(0)$, then  $d_{A,V^*}^{-1}(0)\perp d_{A,V^*}(\B)$ in the sense of G. Birkhoff \cite{BD}), hence $d_{A,V^*}$ has finite ascent \cite [Proposition 2.6]{DDK}.
\begin{thm}
	\label{thm02}  If $d_{A,V^*}^{-1}(0)\subseteq  d_{A^*,V}^{-1}(0)$, then $\a(d_{A,V^*})\leq 1$.
\end{thm} 

\
The following result from \cite{RGD} will be used in some of our argument below. \begin{thm}\label{thm0} If $A, B\in\B$, then the following statements are pairwise equivalent.
 	\vskip3pt\noindent (i) $\rm{ran}(A)\subseteq\rm{ran}(B)$.
 	\vskip3pt\noindent (ii) There is a $\mu\geq 0$ such that $AA^*\leq \mu^2BB^*$.
 	\vskip3pt\noindent (iii) There is an operator $C\in\B$ such that $A=BC$.
 	\vskip4pt\noindent Furthermore, if these conditions hold, then the operator $C$ may be chosen so that (a) $||C||^2=\rm{inf}\{\lambda:AA^*\leq\lambda BB^*\}$; (b) $\rm{ker}(A)=\rm{ker}(C)$; 
 	(c) $\rm{ran}(C)\subseteq\rm{ker}(B)^{\perp}$.\end{thm} 
We note for future reference that $P_m(S,T)=0$ implies $P_m(S^n,T^n)=0$, i.e., {\em $S\in\B$ left $m$-invertible by $T\in\B$ implies $S^n$ left $m$-invertible by $T^n$}, for all $n\in\N$ \cite{DM}.

\

Our main result considers left $m$-invertible operators $S$ such that both $S$ and its left $m$-inverse $T$ are power bounded to prove that such operators are $A$-isometric for some $A>0$.
\begin{thm}
	\label{thm00} If $S\in\B$ is left $m$-invertible by a power bounded operator $T\in\B$, then the following statements are mutually equivalent.
	\vskip4pt\noindent (i) $S$ is power bounded.
	\vskip4pt\noindent (ii) There exists a positive invertible operator $P\in\B$ and an isometry $V\in\B$ such that $S=P^{-1}VP$.
	\vskip4pt\noindent (iii) There exists a positive invertible operator $A\in\B$ such that $T=A^{-1}S^*A$ is a power bounded left $m$-inverse of $S$.
	\vskip4pt\noindent Furthermore, if either of the statements (i), (ii) and (iii) above holds, and $S^*$ has SVEP at $0$ (or, $S$ has a dense range), then $S$ and $T$ are (respectively) similar to some unitaries  $U_1$ and $U_2$  such that $U_1=PU_2P^{-1}$ for some invertible operator $P$.  
\end{thm}

\begin{demo} $(i)\Longrightarrow (ii)$. Let $P_m(S,T)=0$. The hypothesis $S$ and $T$ are power bounded implies the existence of a scalar $M_1>0$ such that $$\sup_{n\in\N}\{||S^n||, ||T^n||\}\leq M_1.$$ The left $m$-invertibility of $S$ by $T$ implies the left invertibility of $S^n$ by $T^n$ for all $n\in\N$, i.e., $$P_m(S,T)=0\Longrightarrow P_m(S^n,T^n)=0, n\in\N$$ \cite{DM}. Define $Z_n$ by $$Z_n=(-1)^{m+1}\sum_{j=0}^m{(-1)^{m-j}\left(\begin{array}{clcr}m\\j\end{array}\right) T^{nj}S^{n(j-1)}}.$$ Then $$P_m(S^n,T^n)=0\Longleftrightarrow Z_nS^n=I$$ for all $n\in\N$ and this, since $$||Z_n||\leq\{1+\left(\begin{array}{clcr}m\\1\end{array}\right)+ \cdots +\left(\begin{array}{clcr}m\\m-1\end{array}\right)\}M^2_1\leq 2^m M_1^2=M$$ for some scalar $M>0$, implies
		$$||x||=||Z_nS^nx||\leq M||S^nx||$$ for all $x\in\H$. Already $$||S^nx||\leq ||S^n||||x||\leq M_1||x||;$$ hence $${\frac{1}{M}}||x||\leq ||S^nx||\leq M_1||x||$$ for all $x\in\H$ and integers $n\geq 1$. Thus, see for example \cite{KR}, $S$ is similar to an isometry $V_1$. Let $$S=EV_1E^{-1}\Longleftrightarrow V_1=E^{-1}SE$$ for some (invertible operator $E\in\B$ and) isometry $V_1$. Then $$S^*|E^{-1}|^2S=|E^{-1}|^2\Longleftrightarrow S^*P^2S=P^2, \hspace{2mm} P=|E|^{-1},$$ implies (by Theorem \ref{thm0}) the existence of an isometry $V\in\B$ such that $$S^*P=PV^*\Longleftrightarrow S=P^{-1}VP.$$

		\noindent $(ii)\Longrightarrow (iii).$ If $S=P^{-1}VP$, $P$ and $V$ as above, then $S^{*n}P^2S^n=P^2$ for all $n\in\N$. Hence $$\sum_{j=0}^m{(-1)^{m-j}\left(\begin{array}{clcr}m\\j\end{array}\right)(P^{-2}S^*P^2)^jS^j}=\sum_{j=0}^m{(-1)^{m-j}\left(\begin{array}{clcr}m\\j\end{array}\right)I}=0,$$ i.e., $P^{-2}S^*P^2=T$ is a power bounded left $m$-inverse of $S$.
			
		\noindent $(iii)\Longrightarrow (i).$ This is immediate from the fact that $P^{-2}S^*P^2$ is power bounded if and only if $S^*$ is power bounded (which in turn holds if and only if $S$ is power bounded).
		
		\vskip4pt Assume next that $(i)$ is satisfied, and hence that $S=P_1U_1P^{-1}_1$ for some isometry $U_1$ and $P_1>0$. If $S^*$ has SVEP (or, $S$ has a dense range), then the left invertibility of $S$ implies $S$ is invertible \cite{A}, and this in turn implies that the isometry $U_1$ is indeed a unitary. Since $$P_m(S,T)=0\Longleftrightarrow P_m(T^*,S^*)=P_m(S,T)^*=0,$$ and the operators $T^*$ and $S^*$ are power bounded, the equivalence $(i)\Longleftrightarrow (iii)$ implies the existence of a positive invertible operator $P_2$ and an isometry $U^*_2$ such that $T^*=P^{-1}_2U^*_2P_2$. Evidently, $T^*$ is ($m$-left invertible, hence) left invertible. We prove that $T^*$, hence $U_2$, is invertible. Clearly,
		\begin{eqnarray*} & & P_m(T^*,S^*)=0\Longleftrightarrow \sum_{j=0}^m{(-1)^{m-j}\left(\begin{array}{clcr}m\\j\end{array}\right)P^{-1}_1U^{*j}_1P_1T^{*j}}=0\\ &\Longleftrightarrow & \sum_{j=0}^m{(-1)^{m-j}\left(\begin{array}{clcr}m\\j\end{array}\right)U^{*j}_1P_1T^{*j}}=0\\ &\Longleftrightarrow& \sum_{j=0}^m{(-1)^{m-j}\left(\begin{array}{clcr}m\\j\end{array}\right)T^jP_1U^j_1}=0.
		\end{eqnarray*}	
	The operator $T^*$ being power bounded has an upper triangular matrix representation
	$$T^*=\left(\begin{array}{clcr} T_1^* & T_0^*\\0 & T_2^*\end{array}\right)\in B(\H_1\oplus\H_2),$$ where $T_1^*\in C_{0.}$ and $T_2^*\in C_{1.}$ \cite{K}. Clearly, $U_2^*=U_{21}^*\oplus U_{22}^*\in B(\H_c\oplus\H_u)$ for some decomposition $\H=\H_c\oplus\H_u$ of $\H$ such that $U_{21}^*=U_2^*|_{\H_c}$ is the backward unilateral shift and $U_{22}^*=U_2^*|_{\H_u}$ is unitary; let $Q=P_2\in B(\H_c\oplus\H_u, \H_1\oplus\H_2)$ have the matrix representation $Q=\left(\begin{array}{clcr}Q_{11} & Q_{12}\\Q_{12}^* & Q_{22}\end{array}\right)$. Then, since $T^*=P_2^{-1}U_2P_2$,  $$Q_{12}^*T_1^*=U_{22}^*Q_{12}^*,\hspace{2mm}Q_{11}T_0^*+Q_{22}T_2^*=Q_{12}U_{22}^*, \hspace{2mm}Q_{11}T^*_1=U_{12}^*Q_{11}.$$ Since  $Q_{12}^*T_1^*=U_{22}Q_{12}^*$ implies $Q_{12}^*T_1^{*n}=U_{22}^{*n}Q_{12}^*$ for all $n\in\N$, since $U_{22}^*$ is unitary and $T_1^*\in C_{0.}$, \begin{eqnarray*}	
	& & ||Q_{12}^*x||= \lim_{n\longrightarrow\infty}||U_{22}^{*n}Q_{12}^*x||\\ &=& \lim_{n\rightarrow\infty}{||Q_{12}^*T_1^{*n}x||}\leq ||Q_{12}^*||\lim_{n\rightarrow\infty}{||T_{1}^{*n}x||}=0\end{eqnarray*} for all $x\in\H_c$. Thus $Q_{12}=0$, and then $Q_{11}$, $Q_{22}$ are invertible positive operators and $Q_{11}T_0^*+Q_{22}T_2^*=Q_{12}U_{22}^*=0$. Furthermore, this from considering $(Q_{11}\oplus Q_{22})T^*=(U_{12}\oplus U_{22})(Q_{11}\oplus Q_{22})$, $Q_{11}T_0^*=0\Longrightarrow T_0^*=0$. Considering now the equation $Q_{11}T^*_1=U_{12}^*Q_{11}$, we have $$||T_1^nx||=||Q_{11}U_{12}^nQ_{11}^{-1}x||\leq ||Q_{11}||||U_{12}^nQ_{11}^{-1}x||$$ for all $n\in\N$ and $x\in\H_c$. Since $U_{12}\in\C_{0.}$, we conclude that $T_1\in C_{0.}\cap C_{.0}=C_{00}$. Hence $T$ is a power bounded operator which is the direct sum of a $C_{00}$ operator with $T_2=Q_{22}^{-1}U_{22}Q_{22}$ (where $U_{22}$ is unitary and $Q_{22}$ is positive invertible). 
	
	\
	
	Define the operator $A\in B(\H_c\oplus\H_u)$, $X\in B(\H_1\oplus\H_2,\H_c\oplus\H_u)$ and $E\in B(\H,\H_c\oplus\H_u)$ by $$A=T_1\oplus U_{22},\hspace{2mm} X=I\oplus Q_{22}\hspace{2mm}\mbox{and}\hspace{2mm} E=XP_1.$$ Then
	
\begin{eqnarray*} & & \sum_{j=0}^m{(-1)^{m-j}\left(\begin{array}{clcr}m\\j\end{array}\right)T^jP_1U^j_1}=0\\ &\Longleftrightarrow& \sum_{j=0}^m{(-1)^{m-j}\left(\begin{array}{clcr}m\\j\end{array}\right)X^{-1}A^jXP_1U^j_1}=0\\ &\Longleftrightarrow& \sum_{j=0}^m{(-1)^{m-j}\left(\begin{array}{clcr}m\\j\end{array}\right)A^jEU^j_1}=0\\ &\Longleftrightarrow& (L_AR_{U_1}-1)^m(E)=\triangle_{A,U_1}^m(E)=0.
\end{eqnarray*}	
	Since the operator $A=XTX^{-1}$ is a power bounded operator which is the direct sum of a unitary with a $C_{.0}$ operator and the operator $U_1$ is  unitary, it follows from an application of Theorems \ref{thm01} and \ref{thm02} that $$\triangle_{A,U_1}(E)=0.$$ Equivalently,
	$$AEU_1-E=0\Longleftrightarrow TP_1U_1-P_1=0\Longleftrightarrow T=P_1U^*_1P_1^{-1}.$$ This implies $T$ is invertible, hence ( $U_2$ is unitary and)  $$P_1^{-1}U_1P_1=P_2U_2P_2^{-1}\Longleftrightarrow U_1=P_1P_2U_2P_2^{-1}P_1^{-1}.$$ Now define (the invertible operator) $P$ by $P=P_1P_2$ to complete the proof.		\end{demo}

It is immediate from the theorem that if a power bounded operator $S_1\in\B$ is left $m$-invertible by a power bounded operator $S_2^*\in\B$, then there exist invertible operators $A_i>0$ in $\B$ such that $S_i$ is an $(A_i,m)$-isometry;  $i=1,2$. 

\

Recall that an operator $A\in\B$ is said to be {\em supercyclic} if there exists a vector $x\in\H$ such that the projective orbit of $x$ under $A$, $${\cal{O}}_{A}({\rm{span}}\{{x}\})=\{{\alpha} {A^n}{x}\in\H:\alpha\in\C,\hspace{2mm} n\geq 0\},$$ is dense in $\H$. Similarities preserve supercyclicity, and power bounded operators of class $C_{1.}$ can not be supercyclic \cite[Theorem 2.1]{AB}. Since isometries are $C_{1.}$ operators:
\begin{cor}
	\label{cor0} If a power bounded operator $S\in\B$ is left $m$-invertible by a power bounded operator $T\in\B$, then neither $S$ nor $T$ is supercyclic. 
\end{cor}

\
Theorem \ref{thm00} is a generalization of the result: $m$-isometric power bounded operators are isometric \cite{D}. That Theorem \ref{thm00} does indeed imply this result is the content of the following proposition. (We remark here that the argument proving the proposition below differs  radically from the argument used in \cite{D}.)
\begin{pro}
	\label{pro00} Power bounded $m$-isometric operators $S\in\B$ are isometric.
\end{pro}
\begin{demo} If $S$ is $m$-isometric and  power bounded, then (as seen above) $S=P^{-1}VP$ for some isometry $V$ and $P>0$. We prove that $[P,S]=PS-SP=0$  (and this would then imply that $S=V$). Evidently, $S=P^{-1}VP$ implies $S^*P=PV^*$. Decompose $V$ into its completely non-unitary (i.e., unilateral shift) and unitary parts by $$V=V_c\oplus V_u\in B(\H_c\oplus\H_u).$$ The operator $S$ being power bounded has an upper triangular matrix representation $$S=\left(\begin{array}{clcr} S_1 & S_0\\0 & S_2\end{array}\right)\in B(\H_1\oplus\H_2),$$ where $S_1\in C_{0.}$ and $S_2\in C_{1.}$ \cite{K}. Let $P\in B(\H_c\oplus\H_u,\H_1\oplus\H_2)$ have the representation $$P=\left(\begin{array}{clcr} P_1 & P_3\\P_3^* & P_2\end{array}\right).$$ Then $S^*P=PV^*$ implies $$S_1^*P_3=P_3V_u^*\Longleftrightarrow V_uP_3^*=P_3^*S_1\Longrightarrow V_u^nP_3^*=P_3^*S_1^n$$ for all $n\in\N$. Hence $$||P_3^*x||=||V^n_uP_3^*x||=||P_3^*S_1^nx|| \leq  ||P_3^*||||S^n_1x||\longrightarrow 0 \hspace{2mm}\mbox{as}\hspace{2mm}n\longrightarrow\infty$$ for all $x\in\H_1$. Thus $P_3=0$ and $P_1,P_2>0$. Since $S^*P=PV^*$ now implies $S_0^*P_1=0$, we must have $S_0=0$. Hence
	$$S=S_1\oplus S_2,\hspace{2mm} S_2=P_2^{-1}V_uP_2\hspace{2mm}\mbox{and}\hspace{2mm}S_1\\in C_{0.}.$$ Evidently, $S$ is $m$-isometric implies $$(L_{S^*_2}R_{S_2}-I)^m(I)=0 \Longleftrightarrow (L_{V_u^*}R_{V_u}-I)^m(P_2^{-2})=0.$$ Applying Theorems 2.1 and 2.2,  it follows that:\begin{eqnarray*}(L_{S^*_2}R_{S_2}-I)^m(I)=0 &\Longleftrightarrow& (L_{V_u^*}R_{V_u}-I)^m(P_2^{-2})=0\\ &\Longleftrightarrow& (L_{V_u^*}R_{V_u}-I)(P_2^{-2})=0\hspace{2mm}(\mbox{i.e.,}\hspace{2mm}\asc(L_{V_u^*}R_{V_u}-I)\leq 1)\\ &\Longleftrightarrow& [V_u,P^{-2}_2]=0\Longleftrightarrow [V_u,P^{-1}_2]=0\\ &\Longrightarrow& S_2=V_u.
\end{eqnarray*} Conclusion: $S^*$ is the direct sum of a $C_{.0}$ and a unitary operator. Applying Theorem \ref{thm01} to $S^*P=PV^*$, $V$ isometric, it follows that $SP=PV$. Hence $$PS^*P=P^2V^*\Longleftrightarrow [V^*,P^2]=0\Longleftrightarrow [V,P]=0\Longleftrightarrow [S,P]=0,$$ and the proof is complete.	
	\end{demo}
The following corollary is immediate from Proposition \ref{pro00}, since either of the hypotheses $S$ has a dense range and $S^*$ has SVEP at $0$ for an $m$-isometric operator $S$ implies the invertibility of $S$.
\begin{cor}\label{cor00} If a power bounded $m$-isometric operator is such that either $S^*$ has SVEP at $0$ or $S$ has a dense range, then  $S$ is unitary.\end{cor}

M. Ch\={o} et al \cite[Theorem 3.15]{CKL} prove that ``if $S\in\B$ is a  power bounded $(m,C)$-isometric operator such that $P_1(CSC,S^*)$ is normaloid (i.e. its norm equals its spectral radius),	then $S\in (1,C)$-isometric. The following proposition shows that the hypothesis $P_1(CSC,S^*)$ is normaloid is redundant, and that the power boundedness of $S$ is sufficient to guarantee $S\in (1,C)$-isometric.

\begin{pro}
	\label{pro01} Power bounded $(m,C)$-isometric operators are $(1,C)$-isometric.
\end{pro}
\begin{demo} By definition,	\begin{eqnarray*} S\in (m,C)-\mbox{isometric} &\Longleftrightarrow& P_m(CSC,S^*)=\sum_{j=0}^m{(-1)^{m-j}\left(\begin{array}{clcr}m\\j\end{array}\right)S^{*j}CS^jC}=0\\ & \Longleftrightarrow & P_m(S,CS^*C)=\sum_{j=0}^m{(-1)^{m-j}\left(\begin{array}{clcr}m\\j\end{array}\right)CS^{*j}CS^j}=0.\end{eqnarray*} Arguing as in the proof of Proposition \ref{pro00}, it follows from the left $m$-invertibility and the power boundedness of $S$ (consequently, also that of $S^*$ and $CS^*C$) that there exists a decomposition $$S=S_1\oplus S_2\in B(\H_1\oplus\H_2),\hspace{2mm} S_1\in C_{0.}\hspace{2mm}\mbox{and}\hspace{2mm} S_2\in C_{1.},$$ of $S$ and a positive invertible operator $Q=Q_1\oplus Q_2\in B(\H_c\oplus \H_u, \H_1\oplus\H_2)$ such that $$S_1=Q^{-1}_1V_cQ_1\hspace{2mm}\mbox{and}\hspace{2mm} S_2=Q^{-1}_2V_uQ_2$$ for some unilateral shift $V_c\in B(\H_c)$ and unitary $V_u\in B(\H_u)$. Set $$Q_1\oplus Q_2=Q, \hspace{4mm}V_c\oplus V_u=V.$$ Evidently, $$S\in(m,C)-\mbox{isometric} \Longleftrightarrow \sum_{j=0}^m{(-1)^{m-j}\left(\begin{array}{clcr}m\\j\end{array}\right)S^{*j}CS^j}= (L_{S^*}R_S-I)^m(C)=0.$$ We prove that $C:\H_1\oplus \H_2\longrightarrow \H_1\oplus\H_2$ has a decomposition $C=C_{11}\oplus C_{22}$. Let $C:{\H}_1\oplus{\H}_2$ into itself have the matrix representation  $$C=\left(\begin{array}{clcr}C_{11} & C_{12}\\C_{21} & C_{22}\end{array}\right).$$  Since
	\begin{eqnarray*}(L_{S^*}R_S-I)^m(C)=0 &\Longleftrightarrow& Q\{(L_{V^*}R_V-I)^m(Q^{-1}CQ^{-1})\}Q=0\\ &\Longleftrightarrow& (L_{V^*}R_V-I)^m(Q^{-1}CQ^{-1})=0,	
	\end{eqnarray*}\begin{eqnarray*} & & (L_{S^*}R_S-I)^m(C)=0\\ &\Longleftrightarrow& \left(\begin{array}{clcr}(L_{V_c^*}R_{V_c}-I)^m(Q_1^{-1}C_{11}Q_1^{-1}) & (L_{V_c^*}R_{V_u}-I)^m(Q_1^{-1}C_{12}Q_2^{-1})\\(L_{V_u^*}R_{V_c}-I)^m(Q_2^{-1}C_{21}Q_1^{-1}) & (L_{V_u^*}R_{V_u}-I)^m(Q_2^{-1}C_{22}Q_2^{-2})\end{array}\right)=0.
\end{eqnarray*}  Set $$(L_{V_c^*}R_{V_u}-I)^{m-1}(Q_1^{-1}C_{12}Q_2^{-1})=Z_{m-1}.$$ Then, $V_u$ ($V_c^*$) being unitary (resp., $C_{0.}$),\begin{eqnarray*}& & (L_{V_c^*}R_{V_u}-I)^m(Q_1^{-1}C_{12}Q_2^{-1})=0\\ &\Longleftrightarrow& (L_{V_c^*}R_{V_u}-I)(Z_{m-1})=V_c^*Z_{m-1}V_u-Z_{m-1}=0\\ &\Longrightarrow& Z_{m-1}=V^{*n}_cZ_{m-1}V_u^n\Longrightarrow ||Z_{m-1}x||=||V_c^{*n} Z_{m-1}V^n_ux||\\ &\longrightarrow& 0\hspace{2mm}\mbox{as}\hspace{2mm} n\longrightarrow\infty
\end{eqnarray*} for all $x\in\H_1$. Hence $Z_{m-1}=0$. Repeating this argument, considering next $(L_{V_c^*}R_{V_u}-I)(Z_{m-2})=0$, a finite number of times it follows that $$Q_1^{-1}C_{12}Q_2^{-1}=0\Longleftrightarrow C_{12}=0.$$ 
A similar argument applied to $$(L_{V_u^*}R_{V_c}-I)^m(T)=0\Longleftrightarrow (L_{V_c^*}R_{V_u}-I)^m(T^*)=0, \hspace{2mm}T=Q^{-1}_2C_{21}Q^{-1}_1,$$ implies that $$T=Q^{-1}_2C_{21}Q^{-1}_1=0\Longleftrightarrow C_{21}=0.$$
Hence $$C=C_{11}\oplus C_{22}.$$
Conside next the equality $(L_{V_u^*}R_{V_u}-I)^m(Q^{-1}_2C_{22}Q_2^{-1})=0$. Set $Q_2^{-1}C_{22}Q_2^{-1}=H$.  Then $(L_{V_u^*}R_{C_{22}V_uC_{22}}-I)^m(HC_{22})=0$. Since $C_{22}V_uC_{22}$ and $V_u^*$ are unitary and $HC_{22}\in B(\H_2)$, $(L_{V_u^*}R_{C_{22}V_uC_{22}}-I)^m(HC_{22})=0$ if and only if $(L_{V_u^*}R_{C_{22}V_uC_{22}}-I)(HC_{22})=0$. Hence

 \begin{eqnarray*} & &   (L_{V_u^*}R_{V_u}-I)^m(Q^{-1}_2C_{22}Q_2^{-1})=0 \\ &\Longrightarrow&   (L_{V_u^*}R_{V_u}-I)(Q^{-1}_2C_{22}Q_2^{-1})=0\\ &\Longleftrightarrow& V_{u}^*Q^{-1}_2C_{22}Q_2^{-1}V_u=Q^{-1}_2C_{22}Q_2^{-1}\\ &\Longleftrightarrow& (Q_2V^*_uQ_2^{-1})C_{22}(Q_2^{-1}V_uQ_2)C_{22}=I\\ &\Longleftrightarrow& S^*_2C_{22}S_2C_{22}-I=0.
\end{eqnarray*} To complete the proof, we prove next that $S_1^*C_{11}S_1C_{11}-I=0$: This would then imply that $$0=(S_1\oplus S_2)^*(C_{11}\oplus C_{22})(S_1\oplus S_2)(C_{11}\oplus C_{22})-I=S^*CSC-I.$$ Set $$(L_{S_1^*}R_{S_1}-I)^{m-1}(C_{11})=X_{m-1}.$$ Then, since $S_1^*$ is power bounded and $S_1\in C_{0.}$, \begin{eqnarray*} & & (L_{S_1^*}R_{S_1}-I)(X_{m-1})=(L_{S_1^*}R_{S_1}-I)^{m}(C_{11})=0\\ &\Longrightarrow& ||X_{m-1}x||=||S_1^{*n}X_{m-1}S^n_1x||\leq ||S_1^{*n}||||X_{m-1}||||S^n_1x||\\ &\longrightarrow& 0\hspace{2mm}\mbox{as}\hspace{2mm}n\longrightarrow\infty
\end{eqnarray*} for all $x\in\H_1$. Hence $X_{m-1}=0$. Repeating the argument, see above,  it follows eventually that $X_1=S^*_1C_{11}S_1-C_{11}=0$. Hence $S^*_1C_{11}S_1C_{11}-I=0$
	\end{demo} 
\begin{rema} {\em As an immediate consequence of Corollary \ref{cor0}, we remark that (just as for $m$-isometries) the similarity of power bounded $(m,C)$-isometric $\B$ operators implies that such operators can not be supercyclic.}\end{rema}


\vskip10pt \noindent\normalsize\rm B.P. Duggal, Faculty of Mathematics, Visgradska 33, 1800 Ni\v{s}, Serbia.\\
\noindent\normalsize \tt e-mail: bpduggal@yahoo.co.uk

\vskip10pt \noindent\normalsize\rm C.S. Kubrusly, Applied Mathematics Department, Federal University, Rio de Janeiro, RJ, Brazil.\\
\noindent\normalsize \tt e-mail:carloskubrusly@gmail.com
\end{document}